\documentclass[12pt]{amsart}
\usepackage{amscd,amssymb,times,fullpage}
\newtheorem{theorem}{Theorem}
\newtheorem{lemma}[theorem]{Lemma}
\newtheorem{proposition}[theorem]{Proposition}
\newtheorem{corollary}[theorem]{Corollary}
\theoremstyle{definition}
\newtheorem{definition}[theorem]{Definition}

\newtheorem{example}[theorem]{Example}

\theoremstyle{remark}
\newtheorem{remark}[theorem]{Remark}

\newcommand{\R}{\mathbb{ R}}

\newcommand{\GG}{\mathcal{ G}}

\newcommand{\FF}{\mathcal{ F}}

\def\bC{{\mathbb C}}

\def\bR{{\mathbb R}}
\def\bH{{\mathbb H}}

\newcommand\Ker{\operatorname{Ker}}

\catcode`\@=\active 

\begin{document}

\title{The geometry of recursion operators}
\author{G.~Bande}
\address{Dipartimento di Matematica e Informatica, Universit\`a degli Studi di Cagliari, Via Ospedale
72, 09129 Cagliari, Italy}
\email{gbande{\char'100}unica.it}
\author{D.~Kotschick}
\address{Mathematisches Institut, Ludwig-Maximilians-Universit\"at M\"unchen,
Theresienstr.~39, 80333 M\"unchen, Germany}
\email{dieter{\char'100}member.ams.org}


\thanks{This work was begun while the second author was a Visiting Professor at the Universit\`a degli Studi di Cagliari,
supported by I.~N.~d.~A.~M. and was completed during a visit of
the first author to Ludwig-Maximilians-Universit\"at M\"unchen
supported in part by G.~N.~S.~A.~G.~A. and P.~R.~I.~N}
\date{October 10, 2007; MSC 2000 classification: primary 53D35, 53C26, 57R30;
secondary 53C12, 53C15, 58A17}

\begin{abstract}
We study the fields of endomorphisms intertwining pairs of symplectic structures. Using these endomorphisms we prove
an analogue of Moser's theorem for simultaneous isotopies of two families of symplectic forms. We also consider the geometric
structures defined by pairs and triples of symplectic forms for which the squares of the intertwining endomorphisms are
plus or minus the identity. For pairs of forms we recover the notions of symplectic pairs and of holomorphic symplectic
structures. For triples we recover the notion of a hypersymplectic structure, and we also find three new structures that
have not been considered before. One of these is the symplectic formulation of hyper-K\"ahler geometry, which turns out
to be a strict generalization of the usual definition in terms of differential or K\"ahler geometry.
\end{abstract}

\maketitle

\section{Introduction}

It is now well known that many geometric structures, particularly on four-manifolds, can be defined in
terms of pairs of two-forms; see for example Donaldson~\cite{D}. In this paper we study the fields of endomorphisms
intertwining such pairs of forms. This leads to a natural generalization of Moser's theorem~\cite{Moser}
on isotopies of symplectic forms and to a generalization of known geometric structures on four-manifolds
to arbitrary dimensions.

Suppose we are given two non-degenerate $2$-forms $\omega$ and $\eta$ on the same manifold $M$. Then there exists
a unique field of invertible endomorphisms $A$ of the tangent bundle $TM$ defined by the equation
\begin{equation}\label{eq:A}
i_X\omega \ = \ i_{AX}\eta \ .
\end{equation}
The important special case when the two $2$-forms involved are closed, and therefore symplectic, is very interesting both
from the point of view of physics, where it arises in the context of bi-Hamiltonian systems, and from a purely mathematical
viewpoint. In physics the field of endomorphisms $A$ is called a recursion operator, and we shall adopt this terminology here.
We shall study the global geometry and topology of a manifold endowed with two (or more) symplectic forms, which we discuss
using the associated recursion operator $A$. For local considerations in the case when the Nijenhuis tensor of $A$ vanishes
see Turiel~\cite{Turiel}.

In Section~\ref{s:Moser} we shall show that the recursion operator neatly encapsulates the necessary and sufficient condition
for the existence of a simultaneous isotopy of two families of symplectic forms. In Section~\ref{s:pairs} we consider the simplest
examples, where the recursion operator $A$ satisfies $A^2=\pm 1$. We shall find that these most basic cases correspond precisely
to the symplectic pairs studied in~\cite{BK}, and to holomorphic symplectic forms respectively. Our discussion of holomorphic
symplectic structures in terms of recursion operators generalizes the work of Geiges~\cite{Geiges2} on conformal symplectic couples
from dimension four to arbitrary dimensions. In Section~\ref{s:triples} we introduce the four geometries defined by triples of symplectic
forms whose pairwise recursion operators all satisfy $A^2=\pm 1$. Throughout our point of view is that of symplectic geometry,
taking as our geometric data only the symplectic forms and the recursion operators they define. Nevertheless, we
shall see that in two of the four cases the data encoded by the triple of symplectic forms define a pseudo-Riemannian metric
leading to the kind of geometry that is used in supersymmetric string theory; see for example~\cite{BGPPR,Hull}.
One of these cases is that of hypersymplectic structures in the sense of Hitchin~\cite{obscure}, the other one is a symplectic analogue
of hyper-K\"ahler structures. We will show that this symplectic formulation of hyper-K\"ahler geometry is not equivalent to the usual one,
because the symplectic data does not force the associated pseudo-Riemannian metric to be definite. Hyper-K\"ahler geometry
corresponds precisely to the special case in which the natural metric is definite. We shall also discuss briefly the geometries defined by triples
of symplectic forms with recursion operators of square $\pm 1$ which do not have natural pseudo-Riemannian metrics attached to  them.
These structures have more to do with foliations than with differential  geometry.

\section{Simultaneous isotopies of symplectic forms}\label{s:Moser}

Given a family $\omega_t$ of smoothly varying symplectic forms on a compact manifold $M$, with $t\in [0,1]$, Moser~\cite{Moser}
showed that there is an isotopy $\varphi_t$ with $\varphi_t^*\omega_t = \omega_0$ if and only if the cohomology class of $\omega_t$
is independent of $t$. This condition ensures that $\dot\omega_t$ is exact, and every choice of a primitive $\alpha_t$ depending
smoothly on $t$ defines a time-dependent vector field $X_t$ by the equation $i_{X_t}\omega_t=-\alpha_t$. The isotopy $\varphi_t$ is
obtained by integrating $X_t$. Conversely, every isotopy with the property $\varphi_t^*\omega_t = \omega_0$ is generated by a
vector field of this form, as one sees by differentiation.

Suppose now that we have two smoothly varying families of symplectic forms $\omega_t$ and $\eta_t$ on a compact manifold $M$.
(We do not make any assumption on the orientations they induce.) When is there an isotopy $\varphi_t$ with $\varphi_t^*\omega_t=
\omega_0$ and $\varphi_t^*\eta_t=\eta_0$? As the vector fields generating isotopies for a single family are very special, one can not
in general expect that there is a vector field which works for both families simultaneously. Let $A$ be the time-dependent recursion
operator defined by
$$
i_X\omega_t \ = \ i_{AX}\eta_t \ .
$$
If an isotopy $\varphi_t$ makes both $\omega_t$ and $\eta_t$ constant, then it makes $A$ constant in $t$. Therefore, an isotopy
can only exist, if the diffeomorphism type of the recursion operator is constant in $t$. If this is the case, we may as well assume
that $A$ is independent of $t$. Then we have the following isotopy result \`a la Moser:
\begin{theorem}\label{t:Moser}
Let $\omega_t$ and $\eta_t$ with $t\in [0,1]$ be smoothly varying families of symplectic forms on a compact manifold $M$, and assume that
the associated recursion operator $A$ is independent of $t$. Then there exists an isotopy $\varphi_t$ with $\varphi_t^*\omega_t=
\omega_0$ and $\varphi_t^*\eta_t=\eta_0$ if and only if $\dot\omega_t$ and $\dot\eta_t$ are exact, and their primitives can be chosen
in such a way that $\dot\omega_t=d\alpha_t$ and $\dot\eta_t=d\beta_t$ with $\alpha_t  = \beta_t\circ A$.
\end{theorem}
\begin{proof}
Suppose that the desired isotopy exists. Then
$$
0 = \frac{d}{dt}\varphi_t^*\omega_t=\varphi_t^*(\dot\omega_t+L_{X_t}\omega_t)=\varphi_t^*(\dot\omega_t+di_{X_t}\omega_t) \ ,
$$
and thus we may take $\alpha_t=-i_{X_t}\omega_t$ as a primitive of $\dot\omega_t$. Similarly we may take
$\beta_t=-i_{X_t}\eta_t$ as a primitive of $\dot\eta_t$. With these choices we have for any $Y\in TM$:
$$
\alpha_t (Y)=-i_{X_t}\omega_t (Y)=-i_{AX_t}\eta_t (Y)=-i_{X_t}\eta_t (AY)=(\beta_t\circ A)(Y)
$$
because $\eta_t (AX,Y)=\eta_t (X,AY)$. Thus the chosen primitives satisfy $\alpha_t  = \beta_t\circ A$.

Conversely, assume that $\dot\omega_t=d\alpha_t$ and $\dot\eta_t=d\beta_t$ with $\alpha_t  = \beta_t\circ A$.
Define two vector fields $X_t$ and $Y_t$ by $i_{X_t}\omega_t=-\alpha_t$ and $i_{Y_t}\eta_t=-\beta_t$. We claim
that $X_t=Y_t$. For the proof we calculate for an arbitrary $Z\in TM$:
$$
i_{X_t}\omega_t (Z) = -\alpha_t (Z)=-\beta_t (AZ) = i_{Y_t}\eta_t (AZ)=i_{AY_t}\eta_t (Z) = i_{Y_t}\omega_t (Z) \ .
$$
The non-degeneracy of $\omega_t$ now implies that $X_t$ and $Y_t$ agree. Denote by $\varphi_t$ the isotopy they
generate. Then $\varphi_t^*\omega_t=\omega_0$ and $\varphi_t^*\eta_t=\eta_0$.
\end{proof}
Theorem~\ref{t:Moser} may look rather ad hoc at first sight, as the geometric meaning of the conditions that $A$ be
independent of $t$ and that it intertwine the primitives $\alpha_t$ and $\beta_t$ is not at all obvious. It may also
not be clear that there are families in which the diffeomorphism type of the recursion operator does change. Nevertheless,
we maintain that this is the natural formulation of the criterion for the existence of simultaneous isotopies. In Section~\ref{s:pairs}
we will specialize this result by making assumptions on $A$, and thereby clarify the geometric content of
Theorem~\ref{t:Moser}. For example, when $A^2=Id_{TM}$ but $A\neq \pm Id_{TM}$, we shall see that $\omega$ and $\eta$
are equivalent to a symplectic pair in the sense of~\cite{BK}, and $A$ contains the information about the pair of foliations
induced by the symplectic pair. In this case it is possible that the foliations could vary in a non-diffeomorphic way, see~\cite{BGK},
so the assumption about the diffeomorphism type of $A$ is not vacuous. Moreover, Theorem~\ref{t:Moser} for this case is
equivalent to the stability theorem for symplectic pairs formulated and proved in~\cite{BGK} using the basic cohomology of foliations.

\section{Symplectic pairs and holomorphic symplectic forms}\label{s:pairs}

The recursion operator $A$ is the identity if and only if $\omega$ and $\eta$ agree. It is minus the identity if and only if
$\omega = -\eta$. From now on we exclude these trivial cases, so we always assume $A\neq\pm Id$.

\subsection{Symplectic pairs}\label{ss:pairs}
Consider first the case $A^2=Id$, but $A\neq\pm Id$. Then the eigenvalues of $A$ are $\pm 1$, and
$$
X \ = \ \frac{1}{2}(X+AX) \ + \ \frac{1}{2}(X-AX)
$$
is the unique decomposition of an arbitrary tangent vector $X$ into a sum of eigenvectors of $A$. Thus the eigenspaces of $A$
give a splitting $TM=D_+\oplus D_-$.
\begin{lemma}
The eigenspaces $D_{\pm}$ for the eigenvalues $\pm 1$ are precisely the kernels of $\Omega^{\mp}=\omega\mp\eta$.
\end{lemma}
\begin{proof}
Let $X$ be an arbitrary tangent vector. Then
$$
i_X\Omega^{\mp} = i_X\omega\mp i_X\eta = i_{AX}\eta \mp i_X\eta = i_{AX\mp X}\eta \ .
$$
As $\eta$ is non-degenerate, the condition $i_X\Omega^{\mp}=0$ is equivalent to $AX=\pm X$.
\end{proof}
The dimensions of the kernels of $\Omega^{\mp}$ are semi-continuous, in that each can only increase on a closed subvariety.
However, the Lemma shows that if the dimension of the kernel of one of the two forms $\Omega^{\mp}$ jumps up, then the
dimension of the kernel of the other one has to decrease. Therefore, the dimensions of the kernels are actually constant on a
connected manifold $M$, so that the forms $\Omega^{\mp}$ have constant ranks. Moreover, as the $\Omega^{\mp}$ are closed,
their kernel distributions are integrable. Thus, the forms $\Omega^{\mp}$ are a symplectic pair in the sense of~\cite{BK}.

Conversely, suppose that we have a symplectic pair $\Omega^{\pm}$ on $M$, that is a pair of closed $2$-forms of constant
ranks, whose kernel foliations $\FF$ and $\GG$ are complementary. Then $\omega = \frac{1}{2}(\Omega^++\Omega^-)$ and 
$\eta= \frac{1}{2}(\Omega^+-\Omega^-)$ are symplectic forms, and the corresponding recursion operator is $A=Id_{T\GG}-Id_{T\FF}$.
Thus $A^2=Id_{TM}$. We have proved:
\begin{theorem}
Two symplectic forms $\omega$ and $\eta$ on a connected manifold $M$ whose recursion operator $A$ satisfies $A^2=Id$
and $A\neq\pm Id$ give rise to a symplectic pair $\Omega^{\pm}$, and every symplectic pair $\Omega^{\pm}$ arises in this way.
\end{theorem}

\begin{remark}
The condition $A^2=Id$ implies that the Nijenhuis tensor of $A$ vanishes identically. Therefore, in this case, $\omega$ and
$\eta$ are compatible in the sense of Poisson geometry.
\end{remark}

The following stability result was proved in~\cite{BGK}:
\begin{theorem}\label{t:BGK}
    Let $\Omega^{\pm}_t$ be a smooth family of symplectic pairs on a closed smooth manifold $M$, such that the kernel
    foliations $\FF=\Ker(\Omega^+_{t})$ and $\GG=\Ker(\Omega^-_{t})$ are independent of $t\in [0,1]$. Then there exists
    an isotopy $\varphi_{t}$ with $\varphi_{t}^{*}\Omega^{\pm}_{t}=\Omega^{\pm}_{0}$ if and only if the basic cohomology
    classes $[\Omega^+_{t}]\in H_{b}^{2}(\FF)$ and $[\Omega^-_{t}]\in H_{b}^{2}(\GG)$ are constant.
    \end{theorem}
We now want to explain the equivalence between this result and Theorem~\ref{t:Moser} in the case when $A^2=Id$.
Consider the symplectic forms $\omega = \frac{1}{2}(\Omega^++\Omega^-)$ and $\eta= \frac{1}{2}(\Omega^+-\Omega^-)$,
and the corresponding recursion operator $A=Id_{T\GG}-Id_{T\FF}$. First of all, assuming that $A$ is independent of $t$ is
the same thing as assuming that its eigenfoliations $\FF$ and $\GG$ are independent of $t$. Example~3.1 of~\cite{BGK}
shows that there are smooth families of symplectic pairs, for which the diffeomorphism type of the foliations is not constant.
In such examples, the diffeomorphism type of the recursion operator is not constant. In the two theorems we assume that
$A$, equivalently $\FF$ and $\GG$, are independent of $t$. The constancy of the cohomology classes of $\omega$ and
$\eta$ is obviously equivalent to the constancy of the cohomology classes of $\Omega^{\pm}_t$, as long as we use de
Rham cohomology in both cases. Now the conditions $\dot\omega_t=d\alpha_t$ and $\dot\eta_t=d\beta_t$ with
$\alpha_t  = \beta_t\circ A$ in Theorem~\ref{t:Moser} are equivalent to the conditions $\dot\Omega^{\pm}_t=d\gamma^{\pm}_t$
with $\gamma^{\pm}_t$ in the ideal of the kernel foliation of $\Omega^{\pm}_t$. This means that the cohomology class
$[\Omega^{\pm}_t]$ is in fact constant in the cohomology of the ideal of the kernel foliation. As explained in~\cite{BGK},
this in turn is equivalent to the constancy of $[\Omega^{\pm}_t]$ in the basic cohomology of the kernel foliation.

\subsection{Holomorphic symplectic structures}\label{s:couples}

Throughout this subsection we assume that we have two symplectic forms $\omega$ and $\eta$ on a manifold $M$ of dimension
$2n$, such that the recursion operator defined by $i_X\omega = i_{AX}\eta$ satisfies $A^2=-Id_{TM}$. This implies
$i_{AX}\omega=-i_X\eta$.

We shall prove the following:
\begin{theorem}\label{t:hol}
If the recursion operator $A$ satisfies $A^2=-Id_{TM}$, then it defines an integrable complex structure with a holomorphic
symplectic form whose real and imaginary parts are $\omega$ and $\eta$. Every holomorphic symplectic form arises in this way.
\end{theorem}
\begin{proof}
In this case $A$ defines an almost complex structure on $M$. We extend $A$ complex linearly to the complexified tangent bundle
$T_{\bC}M=TM\otimes_{\bR}\bC$. The eigenvalues of $A$ are $\pm i$, and
$$
X \ = \ \frac{1}{2}(X-iAX) \ + \ \frac{1}{2}(X+iAX)
$$
is the unique decomposition of a complex tangent vector $X$ into a sum of eigenvectors of $A$. As usual, the eigenspaces of $A$
give a splitting $T_{\bC}M=T^{1,0}\oplus T^{0,1}$, where $T^{1,0}$ is the $+i$ eigenspace, and $T^{0,1}$ is the $-i$ eigenspace. The two
are complex conjugates of each other.
\begin{lemma}
The eigenspaces $T^{0,1}$ and $T^{1,0}$ are precisely the kernels of $\Omega=\omega+i\eta$ and of its complex conjugate
$\bar\Omega=\omega-i\eta$.
\end{lemma}
\begin{proof}
It suffices to prove the statement for the $-i$ eigenspace $T^{0,1}$. The other case then follows by complex conjugation.

Let $X=u+iv$ be a complex tangent vector. Then
$$
i_X\Omega = i_u\omega - i_v\eta + i(i_u\eta+i_v\omega) \ .
$$
The real part of the equation $i_X\Omega=0$ is equivalent to its imaginary part, and each is equivalent to $Au=v$, which is
obviously equivalent to $X\in T^{0,1}$.
\end{proof}

Now we want to see that the almost complex structure $A$ is in fact integrable. By the Newlander--Nirenberg theorem it suffices
to check that one, and hence both, eigendistributions of $A$ are closed under commutation. To do this, suppose $X$ and $Y$
are complex vector fields in $T^{1,0}$, so that $AX=iX$, $AY=iY$. Then, extending $L_X=i_X\circ d + d\circ i_X$ complex linearly to
complex tangent vectors, and using that $\omega$ and $\eta$ are closed, we find
$$
i_{A[X,Y]}\eta=i_{[X,Y]}\omega=L_Xi_Y\omega - i_YL_X\omega=L_Xi_{AY}\eta-i_YL_{AX}\eta = i(L_Xi_{Y}\eta-i_YL_{X}\eta)=i_{i[X,Y]}\eta \ .
$$
The non-degeneracy of $\eta$ now implies that $A[X,Y]=i[X,Y]$, so that in $T^{1,0}$ is closed under commutation.

Thus, we have seen that two symplectic forms $\omega$ and $\eta$ whose recursion operator satisfies $A^2=-Id$ give rise to an
integrable complex structure, for which $T^{0,1}$ is precisely the kernel of $\Omega = \omega +i\eta$. Thus $\Omega$ is a closed
form of type $(2,0)$ and rank $n$, where $n$ is the complex dimension of $M$.

Conversely, if a manifold is complex and carries a holomorphic symplectic form, then the real and imaginary parts of this form
are real symplectic forms whose recursion operator is just the complex structure.

This completes the proof of Theorem~\ref{t:hol}.
\end{proof}

\begin{remark}
The above proof of the integrability of the almost complex structure defined by the recursion operator is the same as that
in Lemma~(6.8) of Hitchin's paper~\cite{HYM}, or in Lemma~(4.1) of~\cite{AH}. However, unlike those references, we do
not assume the symplectic forms to be compatible with any metric. We shall return to a discussion of this in Subsection~\ref{ss:final}.
\end{remark}

The interpretation of Theorem~\ref{t:Moser} in the case of holomorphic symplectic structures is clear. It says that a family of
holomorphic symplectic structures can be made constant by an isotopy if and only if the complex structure $A$ is, up to
diffeomorphism, independent of $t$, the holomorphic symplectic form $\Omega_t$ has constant cohomology class, and the
primitive of $\dot\Omega_t$ can be taken to be a holomorphic form of type $(1,0)$ with respect to the fixed complex structure.

Note that for a manifold with a holomorphic symplectic form the complex dimension $n$ is even, that $\Omega^{n/2}$ is nowhere
zero, and that $\Omega^{(n/2)+1}$ is identically zero. If $n=2$, the latter condition becomes $\Omega^{2}=0$, whose real and
imaginary parts lead to $\omega\wedge\omega=\eta\wedge\eta$ and $\omega\wedge\eta=0$. Thus $\omega$ and $\eta$ precisely
form a conformal symplectic couple in the sense of Geiges~\cite{Geiges2}. A possibly non-conformal couple is a pair of symplectic
forms inducing the same orientation and satisfying the condition $\omega\wedge\eta=0$. For such a couple the recursion operator
may be more complicated, and does not necessarily define a complex structure. However, no closed manifold admitting a
non-conformal couple is known, other than the holomorphic symplectic four-manifolds, which are the $K3$ surfaces, the four-torus,
and the primary Kodaira surfaces. Donaldson~\cite{D} has outlined a strategy that might be applied to prove that every
four-manifold with a symplectic couple also has a conformal one, and is therefore holomorphic symplectic. Note that it follows from
recent work of Li~\cite{Li} that any four-manifold that has a symplectic structure with vanishing first Chern class (as is clearly the case
for the symplectic couples, conformal or not), must have the Betti numbers and intersection form of a holomorphic symplectic four-manifold.

\section{Triples of symplectic forms}\label{s:triples}

We now want to discuss the geometries defined by a triple of symplectic forms $\omega_1$, $\omega_2$, $\omega_3$
whose recursion operators $A_i$ defined by
\begin{equation}
i_X\omega_i = i_{A_{i+2}X}\omega_{i+1} \ ,
\end{equation}
satisfy $A_i^2=\pm Id$ and $A_i\neq\pm Id$. Here and in the sequel all indices are taken modulo $3$. Note that by the
definition all cyclic compositions $A_{i+2} \circ A_{i+1}\circ A_i = Id$.

Depending on how many of the squares of the $A_i$ are $-Id$ and how many are $+Id$, there are four different cases to consider.
We shall see that in the two cases when there is an odd number of $A_i$ with square $-Id$ there are natural pseudo-Riemannian
metrics defined by the triple of two-forms. When exactly one $A_i$ has square $-Id$, we recover the known concept of a
hypersymplectic structure. When all three $A_i$ have square $-Id$, we find a new geometry consisting of a hypercomplex
structure for which all complex structures admit holomorphic symplectic forms. Examples for this new geometric structure, which
we call a hyperholomorphic symplectic structure, are provided by hyper-K\"ahler structures. The latter are precisely those
hyperholomorphic symplectic structures for which the natural pseudo-Riemannian metric is in fact Riemannian. We shall see that
there are non-Riemannian examples as well. In the cases where the number of $A_i$ with square $-Id$ is even there are no natural
metrics, and those structures are rather more flexible than the metric ones.

\subsection{Hyperholomorphic symplectic structures}\label{ss:hyper}

Recall that a hypercomplex structure on a manifold is a triple of integrable complex structures satisfying the quaternion relations;
see for example~\cite{J,S}.

Our first structure given by a triple of symplectic forms is:
\begin{definition}
A triple of symplectic forms $\omega_i$ whose pairwise recursion operators satisfy $A_i^2=-Id$ for all $i=1,2,3$ is called a
hyperholomorphic symplectic structure.
\end{definition}
In this case $A_{i+2} \circ A_{i+1}\circ A_i = Id$ implies that the $A_i$ anti-commute and satisfy the quaternion relations.
By Theorem~\ref{t:hol} each $A_i$ is an integrable complex structure, and so the $A_i$ together form a hypercomplex structure.
Furthermore, each $A_i$ admits a holomorphic symplectic form, justifying the name hyperholomorphic symplectic structure for such
a triple\footnote{This is different from the hypersymplectic structures discussed in~\ref{ss:hypersymplectic} below.}.

There are now many examples of hypercomplex structures, including many on compact manifolds that are not even cohomologically
symplectic; see Sections~7.5 and 7.6 of~\cite{J}, and the references given there. Therefore, hyperholomorphic symplectic structures
are much more restrictive than hypercomplex ones, but, as the following example shows, every hypercomplex structure on $M$ does give rise
to a natural hyperholomorphic symplectic structure on $T^*M$.
\begin{example}
Let $M$ be a manifold with an integrable complex structure $J$. Then lifting $J$ to $T^*M$, the total space of the cotangent
bundle is also a complex manifold. It is also holomorphic symplectic, because if $\omega$ is the exact symplectic form given
by the exterior derivative of the Liouville one-form, then $\Omega(X,Y)=\omega(X,Y)+i\omega(JX,Y)$ is holomorphic symplectic
for the lifted $J$. If $M$ has a hypercomplex structure, then the lifts of the three complex structures to $T^*M$ still
satisfy the quaternion relations, and are the recursion operators for the triple of symplectic forms given by the imaginary parts
of the three holomorphic symplectic forms.
\end{example}

Now we show that hyperholomorphic symplectic structures have natural metrics associated with them.
\begin{proposition}\label{p:hyper}
Let $M$ be a manifold with a hyperholomorphic symplectic structure. Then the bilinear form on $TM$ defined by
$$
g(X,Y)=\omega_i(X,A_iY)
$$
is independent of $i=1,2,3$. It is non-degenerate and symmetric, and invariant under all $A_i$.
\end{proposition}
\begin{proof}
We first prove independence of $i$ as follows:
$$
\omega_i(X,A_iY)= \omega_i(X,A_{i+1}A_{i+2}Y)=\omega_{i+2}(X,A_{i+2}Y)=\ldots=\omega_{i+1}(X,A_{i+1}Y) \ .
$$

Note that $g$ is non-degenerate because $A_i$ is invertible and  $\omega_i$ is non-degenerate.

We prove invariance under the $A_i$ using independence of $i$:
$$
g(A_iX,A_iY)=\omega_{i+1}(A_iX,A_{i+1}A_iY)=
-\omega_{i+1}(A_iX,A_{i+2}Y)=\omega_{i+2}(X,A_{i+2}Y)=g(X,Y) \ .
$$
Finally we prove symmetry using the invariance under $A_i$:
$$
g(Y,X)=\omega_i(Y,A_iX)=-\omega_i(A_iX,Y)=\omega_i(A_iX,A_i^2Y)=g(A_iX,A_iY)=g(X,Y) \ .
$$
\end{proof}
The Proposition shows that $g$ is a pseudo-Riemannian metric compatible with the symplectic forms $\omega_i$.
As it is symmetric and non-degenerate, there must be tangent vectors $X$ with $g(X,X)\neq 0$. Take such a vector
$X$ and consider also $A_1X$, $A_2X$ and $A_3X$. By invariance of $g$ we have $g(A_iX,A_iX)=g(X,X)$, and
by the definition of $g$ and the skew-symmetry of $\omega_i$, the $A_iX$ are $g$-orthogonal to  each other and to $X$.
Replacing $g$ by its negative if necessary, we find
\begin{corollary}\label{c:G}
Every hyperholomorphic symplectic structure in complex dimension two is hyper-K\"ahler.
\end{corollary}
\begin{proof}
Indeed, the pseudo-Riemannian metric $g$ is a definite K\"ahler metric compatible with the underlying hypercomplex
structure, whose K\"ahler forms with respect to $A_i$ are the $\omega_i$ (up to sign).
\end{proof}
\begin{remark}
Interpreting a hyperholomorphic symplectic structure in complex dimension two as a conformal symplectic triple in the
sense of  Geiges~\cite{Geiges2}, Corollary~\ref{c:G} is equivalent to Theorem~2.8 of~\cite{Geiges2}.
\end{remark}

\begin{remark}
For any hyperholomorphic symplectic structure the symplectic forms $\omega_i$ and the complex structures $A_i$
are parallel with respect to the Levi-Civit\`a connection of the pseudo-K\"ahler metric $g$. In particular the Obata
connection of the underlying hypercomplex structure, which is the unique torsion-free connection for which the $A_i$
are parallel, must equal the Levi-Civit\`a connection of $g$.
\end{remark}

In higher dimensions hyper-K\"ahler structures provide examples of hyperholomorphic symplectic structures for
which the natural pseudo-Riemannian metric $g$ is definite. However, there are many other examples, even on
manifolds that do not support any K\"ahler structure, so that Corollary~\ref{c:G} does not generalize to higher dimensions,
as shown by the following result.
\begin{theorem}\label{t:nonK}
In every even complex dimension $\geq 4$ there exist hyperholomorphic symplectic structures on closed
manifolds that do not support any K\"ahler structure.
\end{theorem}
\begin{proof}
Our examples will be nilmanifolds. Among such manifolds, it is known that only tori admit K\"ahler structures;
see Benson and Gordon~\cite{BG}. Therefore it is enough to produce a nilmanifold of real dimension $8$ that
is not a torus but admits a hyperholomorphic symplectic structure. Then we can take products with $T^4$ to
prove the result in all dimensions.

Our eight-dimensional example comes from the work of Dotti and  Fino~\cite{DF}, who found a non-Abelian two-step nilpotent
Lie algebra which is both hypercomplex and symplectic. What is new here, is that we write down three invariant symplectic forms
such that the recursion operators are complex structures forming a hyperholomorphic symplectic structure.

Consider the real Lie algebra $\mathfrak{g}$ spanned by $8$ vectors $e_1,\ldots,e_8$ such that
$$
[e_1,e_3]=-[e_2,e_4]=e_7, \ \ \ \ \ \ [e_1,e_4]=[e_2,e_3]=e_8 \ ,
$$
and all other commutators vanish. Clearly this is nilpotent. As the structure constants are rational, the corresponding simply
connected nilpotent Lie group $G$ admits cocompact discrete subgroups $\Gamma$, and our example will be $M=G/\Gamma$.
This is a nilmanifold, and is not a torus because $\mathfrak{g}$ is not Abelian.

We can take a framing of $G$ by left-invariant vector fields
corresponding to the $e_i$. Let $e^i$ be the dual left-invariant
one-forms. Then $e^1,\ldots,e^6$ are closed and, by the above
formulae, we have
$$
de^7 = -e^1\wedge e^3+e^2\wedge e^4 , \ \ \ \ \ \ de^8 = -e^1\wedge e^4 - e^2\wedge e^3 \ .
$$
Now we consider the following left-invariant two-forms:
$$
\omega_1 = e^8\wedge e^1 + e^7\wedge e^2 - e^6\wedge e^3 + e^5\wedge e^4 \ ,
$$
$$
\omega_2 = e^8\wedge e^2 - e^7\wedge e^1 + e^6\wedge e^4 + e^5\wedge e^3 \ ,
$$
$$
\omega_3 = e^8\wedge e^3 + e^7\wedge e^4 + e^6\wedge e^1 - e^5\wedge e^2 \ .
$$
These forms are clearly non-degenerate, and by substituting from the formulae  for $de^7$ and $de^8$ we see that they are closed.
A direct calculation shows that the recursion operators are almost complex structures.
Thus, by our previous discussion, we have
a left-invariant hyperholomorphic symplectic structure, which descends to $G/\Gamma$.
\end{proof}
In the $8$-dimensional example used in the proof the pseudo-Riemannian metric $g$ has signature $(4,4)$.
By taking products of this and of hyper-K\"ahler examples, we can realize all possible signatures of the form $(4k,4l)$
with $k+l\geq 2$ as signatures of hyperholomorphic symplectic structures.
Other examples can be constructed using the nilpotent Lie algebras also used in~\cite{FPPS}.

\subsection{Hypersymplectic structures}\label{ss:hypersymplectic}

Next we consider a triple of symplectic forms such that two recursion operators have square the identity, and one has square
minus the identity. After renumbering we may assume $A_1^2=-Id$ and $A_2^2=A_3^2=Id$. Then the cyclic relations
$A_{i+2} \circ A_{i+1}\circ A_i = Id$ show that
the $A_i$ anti-commute and $A_2A_1=A_3$. It follows that $A_i\neq\pm Id$, so the trivial cases are excluded automatically.

We have the following result analogous to Proposition~\ref{p:hyper}:
\begin{proposition}\label{p:hypersymplectic}
Let $M$ be a manifold with three symplectic forms whose recursion operators satisfy $A_1^2=-Id$ and $A_2^2=A_3^2=Id$.
Then
$$
\omega_1(X,A_1Y)=-\omega_2(X,A_2Y)=-\omega_3(X,A_3Y) \ ,
$$
and these expressions define a bilinear form $g(X,Y)$ on $TM$. It is non-degenerate and symmetric, invariant under $A_1$, and
satisfies $g(A_iX,A_iY)=-g(X,Y)$ for $i=2,3$.
\end{proposition}
We omit the proof as it is literally the same as for Proposition~\ref{p:hyper}.

Now in this case if we take a vector $X$ with $g(X,X)\neq 0$, then $g(A_1X,A_1X)=g(X,X)$, and $g(A_2X,A_2X)=g(A_3X,A_3X)=-g(X,X)$,
and the $A_iX$ are $g$-orthogonal to each other and to $X$. Thus, we have a $4$-dimensional subspace on which
$g$ is non-degenerate and has signature $(2,2)$. Looking at the orthogonal complement of this subspace and proceeding
inductively, we see that the metric $g$ has neutral signature.

We can compare this data with the following definition due to Hitchin~\cite{obscure}; see also~\cite{DJS,FPPS}.
\begin{definition}
A hypersymplectic structure on a manifold is a pseudo-Riemannian metric $g$ of neutral signature, together with three
endomorphisms $I$, $S$ and $T$ of the tangent bundle satisfying
$$
I^2=-Id \ , \ \ \ S^2=T^2=Id \ , \ \ \ IS=-SI=T \ ,
$$
$$
g(IX,IY)=g(X,Y) \ , \ \ \ g(SX,SY)=-g(X,Y) \ , \ \ \ g(TX,TY)=-g(X,Y) \ ,
$$
and such that the following three two-forms are closed:
$$
\omega_I(X,Y)=g(IX,Y) \ , \ \ \ \omega_S(X,Y)=g(SX,Y) \ , \ \ \ \omega_T(X,Y)=g(TX,Y) \ .
$$
\end{definition}
Given a hypersymplectic structure in this sense, the recursion operators intertwining the three symplectic forms are,
up to sign, precisely the endomorphisms $I$, $S$ and $T$. Conversely, given three symplectic forms for which one of the
pairwise recursion operators has square $-Id$ and the other two have square the identity, Proposition~\ref{p:hypersymplectic}
shows that we can recover a uniquely defined hypersymplectic structure. Thus we have proved:
\begin{corollary}\label{c:equiv}
A hypersymplectic structure is equivalent to a unique triple of symplectic forms for which two of the recursion operators have
square the identity, and one has square minus the identity.
\end{corollary}

In real dimension $4$ we have the following classification of closed hypersymplectic manifolds, which one can think of
as a hypersymplectic analogue of Corollary~\ref{c:G}.
\begin{proposition}
[cf.~\cite{Kamada}]
A closed oriented four-manifold admits a hypersymplectic structure if and only it is $T^4$ or a nilmanifold for
$Nil^3\times\R$.
\end{proposition}
\begin{proof}
A closed oriented four-manifold with a hypersymplectic structure is holomorphic symplectic, and so by a result of Kodaira is $T^4$, a
primary Kodaira surface, or a $K3$ surface; see~\cite{BPV}. Clearly $T^4$ inherits the standard hypersymplectic structure of $\R^4$.

By a result of  Wall~\cite{W2}, primary Kodaira surfaces are precisely the nilmanifolds of $Nil^3\times\R$. Recall from~\cite{W2}
or~\cite{BK} that $Nil^3\times\R$ has a framing by left-invariant one-forms $\alpha_1,\ldots,\alpha_4$ with
$d\alpha_3=\alpha_1\wedge\alpha_2$, and $\alpha_i$ closed for $i\neq 3$. The left-invariant two-forms
$$
\omega_1 = \alpha_3\wedge\alpha_1 + \alpha_2\wedge\alpha_4
$$
$$
\omega_2 = \alpha_3\wedge\alpha_2 - \alpha_1\wedge\alpha_4
$$
$$
\omega_3 = \alpha_3\wedge\alpha_2 + \alpha_1\wedge\alpha_4
$$
define an invariant hypersymplectic structure that descends to all compact quotients.

A hypersymplectic structure also defines a symplectic pair, and therefore a four-manifold with such a structure is symplectic
for both choices of orientation. But a $K3$ surface endowed with the non-complex orientation can not be symplectic, because
it has vanishing Seiberg--Witten invariants. This follows from the existence of smoothly embedded spheres whose selfintersection
number is positive for the non-complex orientation; cf.~\cite{Ko}.
\end{proof}

High-dimensional examples of hypersymplectic structures on closed manifolds have recently appeared in~\cite{FPPS,AD}.

\subsection{Holomorphic symplectic pairs}

Next we consider the following:
\begin{definition}
A triple of symplectic forms $\omega_i$ is called a holomorphic symplectic pair if two of the pairwise recursion operators
have square $-Id_{TM}$ and one has square $Id_{TM}$, but is not itself $\pm Id_{TM}$.
\end{definition}
After renumbering we may assume $A_1^2=A_2^2=-Id$ and $A_3^2=Id$. Then the cyclic relations
$A_{i+2} \circ A_{i+1}\circ A_i = Id$ show that the $A_i$ commute and $A_2A_1=A_3$.

Now $A_3$ has square the identity, but is not itself plus or minus the identity, and so defines a symplectic pair.
The other two recursion operators, $A_1$ and $A_2$ define integrable complex structures. As they commute with
$A_3$, they preserve its eigenfoliations and restrict as complex structures to the leaves. On the $+1$ eigenfoliation
of $A_3$ the two complex structures are complex conjugates of each other, and on the $-1$ eigenfoliation they agree.
The two complex structures also have holomorphic symplectic forms which restrict to the leaves of the
eigenfoliations of $A_3$. Thus a holomorphic symplectic pair is a symplectic pair whose leaves are not just symplectic,
but are holomorphic symplectic submanifolds. It follows in particular that the real dimensions of the leaves are multiples
of $4$, and the smallest dimension in which this structure can occur is $8$.

Here are some examples.
\begin{example}
Consider the eight-dimensional nilpotent Lie group $G$ from the proof of Theorem~\ref{t:nonK}. The forms
$$
\omega_1 = e^8\wedge e^1 + e^7\wedge e^2 - e^6\wedge e^3 + e^5\wedge e^4
$$
$$
\omega_2 = e^8\wedge e^2 - e^7\wedge e^1 + e^6\wedge e^4 + e^5\wedge e^3
$$
$$
\omega_3 = e^8\wedge e^2 - e^7\wedge e^1 - e^6\wedge e^4 - e^5\wedge e^3
$$
are a left-invariant holomorphic symplectic pair, that descends to all compact quotients.
\end{example}

\begin{example}
Let $M_1$ and $M_2$ be holomorphic symplectic manifolds, and let $(\omega_i,\eta_i)$ be the two symplectic forms
defining the structure on $M_i$. On the product $M_1\times M_2$ the three symplectic forms
$\omega_1+\omega_2$, $\eta_1+\eta_2$ and $-\omega_1+\omega_2$ are a holomorphic symplectic pair. In this case
the foliations are given by the factors of the product.
\end{example}
For the factors in this construction we can use any holomorphic symplectic manifold. This could be hyper-K\"ahler, or
a nilmanifold, or one of the simply connected non-K\"ahler examples constructed by Guan~\cite{Guan}, compare also~\cite{Bog}.

\subsection{Symplectic triples}

We use the name symplectic triple\footnote{This terminology clashes with that of Geiges~\cite{Geiges2}.}
 for a triple of symplectic forms such that $A_i^2= Id$ and $A_i\neq\pm Id$ for all three $i$.
By the discussion in Subsection~\ref{ss:pairs} above this is equivalent  to requiring that for $i\neq j$, the forms $\omega_i$
and $\omega_j$ define a symplectic pair. In other words, $\omega_i\pm\omega_j$ are forms of constant (non-maximal and
non-zero) rank, and each one restricts as a symplectic form to the leaves  of the kernel foliation of the other one. In
dimension four, the two symplectic forms making up a symplectic pair induce opposite orientations. This means in particular
that there can not be any symplectic triples. Starting in dimension $6$, however, symplectic triples exist in abundance;
see~\cite{BK}, particularly Remark~7. Using the examples of symplectic pairs on four-manifolds constructed in~\cite{BK},
one immediately obtains many examples of symplectic triples in higher dimensions by taking products with other symplectic
manifolds. In this way many different topological types can be realized.

\begin{example}
The simplest example is given by considering three closed two-forms $\eta_i$ of constant rank $=2$ on a $6$-manifold, with
the property that $\eta_1\wedge\eta_2\wedge\eta_3$ is nowhere zero. Then we can take $\omega_1=\eta_1+\eta_2+\eta_3$,
$\omega_2=\eta_1+\eta_2-\eta_3$ and $\omega_3=\eta_1-\eta_2-\eta_3$. This works for example by taking a product of three
surfaces, or, more interestingly, by taking a quotient of the polydisk $\bH^2\times\bH^2\times\bH^2$ by an irreducible lattice as
in Subsection~5.2 of~\cite{BK}.
\end{example}

\subsection{Final comments and remarks}\label{ss:final}

We have seen in Corollary~\ref{c:equiv} that the metric definition
of hypersymplectic structures given in~\cite{obscure} is in fact
equivalent to the symplectic definition we have given in terms of
symplectic forms and their recursion operators. In particular, the
fact that the signature of the natural pseudo-K\"ahler metric is
always neutral follows purely from the algebraic relations between
the recursion operators. The analogous result is false for
hyper-K\"ahler structures. If one retains from a hyper-K\"ahler
structure only the triple of symplectic forms together with the
property that the intertwining recursion operators be almost
complex structures, then one finds hyperholomorphic symplectic
structures. As we have seen in Subsection~\ref{ss:hyper}, the
ensuing algebraic identities define a pseudo-Riemannian metric of
a priori unknown signature. Therefore, when viewing hyper-K\"ahler
structures as symplectic objects, as is done for example
in~\cite{obscure}, one either has to allow pseudo-hyper-K\"ahler
structures, or one has to build the definiteness of the metric
into the definition. It is not clear to us how one would do this
in terms of symplectic geometry alone. In this direction, the
discussion on page 172 of~\cite{obscure} is really based on
Lemma~(6.8) of~\cite{HYM}, where definiteness of the compatible
metric is part of the definition.

Finally, for holomorphic symplectic pairs and for symplectic triples there are no natural metrics, definite or otherwise.
We leave it to the interested reader to work out why Proposition~\ref{p:hyper} has no analogue in these cases.

\subsection*{Acknowledgement}
We are grateful to N.~Hitchin for some useful comments and for sending us a copy of~\cite{obscure}.

\bigskip

\bibliographystyle{amsplain}


\end{document}